\documentclass[a4paper,9pt]{article}
\usepackage{amsmath, amsfonts, amssymb}
\usepackage[english, russian]{babel}
\usepackage[cp1251]{inputenc}

\usepackage{graphicx}
\begin{document}
\sloppy
\begin{center}
\textbf{Some relations following from the decomposition formula \\ for one
multidimensional Lauricella hypergeometric function }\\[5pt]
\textbf{Ergashev T.G.\\}
\medskip
{ Institute of Mathematics, Uzbek
Academy of Sciences,  Tashkent, Uzbekistan. \\

{\verb ergashev.tukhtasin@gmail.com }}\\
\end{center}

\begin{quote} 

Fundamental solutions for a class of multidimensional elliptic
equations with several singular coefficients were constructed recently. These fundamental
solutions  are directly connected with  multiple Lauricella hypergeometric
function and  the decomposition formula is required for their
investigation which would express the multivariable hypergeometric
function in terms of products of several simpler hypergeometric
functions involving fewer variables. In this paper, some relations following from the decomposition formula for one
multidimensional Lauricell hypergeometric function are determined.
\end{quote}

We consider the equation

\begin{equation}
\label{eq1}
L_{\alpha} ^{(m,n)} \left( {u} \right): = {\sum\limits_{i = 1}^{m} {u_{x_{i}
x_{i}}} }   + {\sum\limits_{j = 1}^{n} {{\frac{{2\alpha _{j}}} {{x_{j}
}}}u_{x_{j}}} }   = 0
\end{equation}

\noindent
in the region $R_{m}^{n +}  = {\left\{ {x:x_{1} > 0,x_{2} > 0,...,x_{n} > 0}
\right\}},$ where $x = \left( {x_{1} ,...,x_{m}}  \right), \, m \ge 2,\\ 0 \le
n \le m;\alpha = \left( {\alpha _{1} ,...,\alpha _{n}}  \right),\alpha
_{j} $ are real numbers with $0 < 2\alpha _{j} < 1,j = \overline {1,n} .$

\bigskip

Fundamental solutions of the equation were constructed recently \cite{Erg}. In fact,
the fundamental solutions of the equation (\ref{eq1}) can be expressed in terms of
Lauricella's hypergeometric function in $n\,$variables, that is, the
Lauricella multivariable hypergeometric function
\[
F_{A}^{(n)} \left( {a,b_{1} ,...,b_{n} ;c_{1} ,...,c_{n} ;z_{1} ,...,z_{n}}
\right)
\]
\noindent
defined by
\[
F_{A}^{\left( {n} \right)} \left( {a,b_{1} ,...,b_{n};c_{1} ,...,c_{n};z_{1} ,...,z_{n}} \right) = F_{A}^{(n)} {\left[
{{\begin{array}{*{20}c}
 {a,b_{1} ,...,b_{n}}  \hfill \\
 {c_{1} ,...,c_{n}}  \hfill \\
\end{array}} z_{1} ,...,z_{n}}  \right]}
\]
\begin{equation}
\label{eq2}
 = {\sum\limits_{p_{1} ,...,p_{n} = 0}^{\infty}  {{\frac{{\left( {a}
\right)_{p_{1} + ... + p_{n}}  \left( {b_{1}}  \right)_{p_{1}}  ...\left(
{b_{n}}  \right)_{p_{n}}} } {{\left( {c_{1}}  \right)_{p_{1}}  ...\left(
{c_{n}}  \right)_{p_{n}}} } }}} {\frac{{z_{1}^{p_{1}}} } {{p_{1}
!}}}...{\frac{{z_{n}^{p_{n}}} } {{p_{n} !}}},
\quad
{\sum\limits_{i = 1}^{n} {{\left| {z_{i}}  \right|} < 1,}}
\end{equation}
\noindent
where $c_{i} \ne 0, - 1, - 2,...,\,\,i = \overline {1,n} $ and $\left(
{\kappa}  \right)_{\nu}  $ denotes the general Pochhammer symbol or the
shifted factorial, since
\[
\left( {1} \right)_{m} = m!
\left( {m \in {\rm N}_{0} : = {\rm N} \cup \{0\};\,\,{\rm N} = {\left\{
{1,2,3,...} \right\}}} \right),
\]
\noindent
which is defined in terms of the familiar Gamma function, by
\[
\left( {\kappa}  \right)_{0} = 1,
\left( {\kappa}  \right)_{\nu}  = {\frac{{\Gamma \left( {\kappa + \nu}
\right)}}{{\Gamma \left( {\kappa}  \right)}}},
\nu \in {\rm N}.
\]

We thus obtain the following fundamental solutions:
\begin{equation}
\label{eq3}
\begin{array}{l}
 q_{k} \left( {x,\xi}  \right) = \gamma _{k} {\prod\limits_{i = 1}^{k}
{\left( {x_{i} \xi _{i}}  \right)^{1 - 2\alpha _{i}}  \cdot}}  r^{ - 2\bar
{\alpha} _{k}}
 F_{A}^{\left( {n} \right)} {\left[ {{\begin{array}{*{20}c}
 {\bar {\alpha} _{k} ,1 - \alpha _{1} ,...,1 - \alpha _{k} ,\alpha _{k + 1}
,...,\alpha _{n} ;} \hfill \\
 {2 - 2\alpha _{1} ,...,2 - 2\alpha _{k} ,2\alpha _{k + 1} ,...,2\alpha _{n}
;} \hfill \\
\end{array}} \sigma _{1} ,...,\sigma _{n}}  \right]}, \\
 \end{array}
\end{equation}

\noindent
where
\[
\bar {\alpha} _{k} = {\frac{{m}}{{2}}} + k - 1 - {\sum\limits_{i = 1}^{k}
{\alpha _{i}}}   + {\sum\limits_{i = k + 1}^{n} {\alpha _{i}}}  ,
\]
\[
\gamma _{k} = 2^{2\bar {\alpha} _{k} - m}{\frac{{\Gamma \left( {\bar {\alpha
}_{k}}  \right)}}{{\pi ^{m / 2}}}}{\prod\limits_{i = k + 1}^{n}
{{\frac{{\Gamma \left( {\alpha _{i}}  \right)}}{{\Gamma \left( {2\alpha _{i}
} \right)}}}}} {\prod\limits_{j = 1}^{k} {{\frac{{\Gamma \left( {1 - \alpha
_{j}}  \right)}}{{\Gamma \left( {2 - 2\alpha _{j}}  \right)}}}}} ,
\quad
k = \overline {0,n} ;
\]
\[
\sigma _{k} = 1 - {\frac{{r_{k}^{2}}} {{r^{2}}}},
\quad
r^{2} = {\sum\limits_{i = 1}^{m} {\left( {x_{i} - \xi _{i}}  \right)^{2}}} ,
\quad
r_{k}^{2} = \left( {x_{k} + \xi _{k}}  \right)^{2} + {\sum\limits_{i = 1,i
\ne k}^{m} {\left( {x_{i} - \xi _{i}}  \right)^{2}}} ,
k = \overline {1,n}.
\]
Here ${\sum\limits_{i = l + 1}^{l} {}} $ is to be interpreted as zero, when $l = 0$
or $l = n$, and ${\prod\limits_{i = l + 1}^{l} {}} $ is to be interpreted as one, when $l = 0$ or $l = n$.

For a given multiple hypergeometric function, it is useful to fund a
decomposition formula which would express the multivariable hypergeometric
function in terms of products of several simpler hypergeometric functions
involving fewer variables. Burchnall and Chaundy \cite{{BC1},{BC2}} systematically
presented a number of expansion and decomposition formulas for some double
hypergeometric functions in series of simpler hypergeometric functions. For
example, the Appell function
\begin{equation*}
F_{2} \left( {a,b_{1} ,b_{2} ;c_{1} ,c_{2} ;x,y} \right) = {\sum\limits_{i,j
= 0}^{\infty}  {{\frac{{\left( {a} \right)_{i+j} \left( {b_{1}}
\right)_{i} \left( {b_{2}}  \right)_{j}}} {{\left( {c_{1}}  \right)_{i}
\left( {c_{2}}  \right)_{j}}} }{\frac{{x^{i}}}{{i!}}}{\frac{{y^{j}}}{{j!}}}}
}
\end{equation*}
\begin{equation*}
{\left[ {c_{1} ,c_{2} \ne 0,-1,-2,...;\,\,{\left| {x} \right|} + {\left| {y} \right|}
< 1} \right]}
\end{equation*}
\noindent
has the expansion \cite{BC1}
\begin{equation*}
F_{2} \left( {a,b_{1} ,b_{2} ;c_{1} ,c_{2} ;x,y} \right)\]
\[ = {\sum\limits_{i =
0}^{\infty}  {{\frac{{\left( {a} \right)_{i} \left( {b_{1}}  \right)_{i}
\left( {b_{2}}  \right)_{i}}} {{i!\left( {c_{1}}  \right)_{i} \left( {c_{2}
} \right)_{i}}} }x^{i}y^{i}F\left( {a + i,b_{1} + i;c_{1} + i;x}
\right)F\left( {a + i,b_{2} + i;c_{2} + i;y} \right)}} ,
\end{equation*}
where
$$F {\left(a,b;c;z\right)}=F {\left[
{{\begin{array}{*{20}c}
 {a,b} ; \hfill \\
 {c};  \hfill \\
\end{array}} z}  \right]}={\sum\limits_{i =
0}^{\infty} \frac{(a)_i(b)_i}{(c)_ii!}z^i }$$
is Gaussian hypergeometric function \cite{Erd}.

The Burchnall-Chaundy method, which is limited to functions of two variables,
is based on the following mutually inverse symbolic operators \cite{BC1}
\begin{equation}
\label{eq8}
\nabla \left( {h} \right) = {\frac{{\Gamma \left( {h} \right)\Gamma \left(
{{\rm \delta} _{1} + {\rm \delta} _{2} + h} \right)}}{{\Gamma \left( {{\rm
\delta} _{1} + h} \right)\Gamma \left( {{\rm \delta} _{2} + h} \right)}}},
\quad
\Delta \left( {h} \right) = {\frac{{\Gamma \left( {{\rm \delta} _{1} + h}
\right)\Gamma \left( {{\rm \delta} _{2} + h} \right)}}{{\Gamma \left( {h}
\right)\Gamma \left( {{\rm \delta} _{1} + {\rm \delta} _{2} + h}
\right)}}},
\end{equation}
\noindent
where ${\rm \delta} _{1} = x{\displaystyle\frac{{\partial}} {{\partial x}}}$ and ${\rm
\delta} _{2} = y{\displaystyle\frac{{\partial}} {{\partial y}}}$.

In order to generalize the operators $\nabla \left( {h} \right)$ and $\Delta
\left( {h} \right)$, defined in (4), A.Hasanov and H.M.Srivastava \cite{{HS6}, {HS7}}
introduced the operators
\begin{equation}
\label{eq9}
\tilde {\nabla} _{z_{1} ;z_{2} ,...,z_{n}}  \left( {h} \right) =
{\frac{{{\rm \Gamma} \left( {h} \right){\rm \Gamma} \left( {{\rm \delta
}_{1} + ... + {\rm \delta} _{n} + h} \right)}}{{{\rm \Gamma} \left( {{\rm
\delta} _{1} + h} \right){\rm \Gamma} \left( {{\rm \delta} _{2} + ... + {\rm
\delta} _{n} + h} \right)}}},
\end{equation}
\begin{equation}
\label{eq10}
\tilde {\Delta} _{z_{1} ;z_{2} ,...,z_{n}}  \left( {h} \right) =
{\frac{{{\rm \Gamma} \left( {{\rm \delta} _{1} + h} \right){\rm \Gamma
}\left( {{\rm \delta} _{2} + ... + {\rm \delta} _{n} + h} \right)}}{{{\rm
\Gamma} \left( {h} \right){\rm \Gamma} \left( {{\rm \delta} _{1} + ... +
{\rm \delta} _{n} + h} \right)}}},
\end{equation}
\noindent
where ${\rm \delta} _{k} = z_{k} {\displaystyle\frac{{\partial}} {{\partial z_{k}}} }\,,$
with the help of which they managed to find decomposition formulas for a
whole class of hypergeometric functions in several variables. For example,
the hypergeometric Lauricella function $F_{A}^{\left( {n} \right)} $,
defined by formula (\ref{eq2}) has the decomposition formula \cite{HS6}
\begin{equation*}
F_{A}^{(n)} \left( {a,b_{1} ,...,b_{n} ;c_{1} ,...,c_{n} ;z_{1} ,...,z_{n}}
\right)
\end{equation*}
\begin{equation*}
= {\sum\limits_{m_{2} ,...,m_{n} = 0}^{\infty}  {}} {\frac{{\left(
{a} \right)_{m_{2} + ... + m_{n}}  \left( {b_{1}}  \right)_{m_{2} + ... +
m_{n}}  \left( {b_{2}}  \right)_{m_{2}}  ...\left( {b_{n}}  \right)_{m_{n}}
}}{{m_{2} !...m_{n} !\left( {c_{1}}  \right)_{m_{2} + ... + m_{n}}  \left(
{c_{2}}  \right)_{m_{2}}  ...\left( {c_{n}}  \right)_{m_{n}}
}}}z_{1}^{m_{2} + ... + m_{n}}  z_{2}^{m_{2}}  ...z_{n}^{m_{n}}
\end{equation*}
\begin{equation*}
 \cdot F\left( {a + m_{2} + ... + m_{n} ,b_{1} + m_{2} + ... + m_{n} ;c_{1}
+ m_{2} + ... + m_{n} ;z_{1}}  \right)
\end{equation*}
\begin{equation*}
 \cdot F_{A}^{(n - 1)} ( a + m_{2} + ... + m_{n} ,b_{2} + m_{2}
,...,b_{n} + m_{n} ;
\end{equation*}
\begin{equation}
\label{eq11}
c_{2} + m_{2} ,...,c_{n} + m_{n} ;z_{2} ,...,z_{n}
),n \in {\rm N}\backslash {\left\{ {1} \right\}}.
\end{equation}

However, due to the recurrence of formula (7), additional difficulties may
arise in the applications of this expansion. Further study of the properties
of operators (5) and (6) showed that formula (7) can be reduced to a
more convenient form.

\bigskip

\textbf{Lemma 1} \cite{Erg}. The following decomposition formula holds true at $n
\in {\rm N}\backslash {\left\{ {1} \right\}}$
\begin{equation*}
F_{A}^{(n)} \left( {a,b_{1} ,b_{2} ,....,b_{n} ;c_{1} ,c_{2} ,....,c_{n}
;z_{1} ,...,z_{n}}  \right)
\end{equation*}
\begin{equation*}
 = {\sum\limits_{{\mathop {m_{i,j} = 0}\limits_{(2 \le i \le j \le n)}
}}^{\infty}  {{\frac{{(a)_{A(n,n)}}} {{m_{ij} !}}}}} \prod\limits_{k =
1}^{n} [ \frac{{(b_{k} )_{B(k,n)}}}{(c_{k})_{B(k,n)}}
 z_{k}^{B(k,n)} F( a + A(k,n),\end{equation*}
\begin{equation}
\label{eq12} b_{k} + B(k,n);c_{k} + B(k,n);z_{k})] ,
\end{equation}
\noindent
where
\begin{equation}
\label{eq13}
A(k,n) = {\sum\limits_{i = 2}^{k + 1} {{\sum\limits_{j = i}^{n} {m_{i,j}}
}}} , \,\,
B(k,n) = {\sum\limits_{i = 2}^{k} {m_{i,k} +}}  {\sum\limits_{i = k + 1}^{n}
{m_{k + 1,i}}}.
\end{equation}

The formula (8) is proved by the method mathematical induction \cite{Erg}.

It should be noted here that the sum ${\sum\limits_{k = 1}^{n} {B(k,n)}} $
has the parity property, which plays an important role in the calculation of
the some values of hypergeometric functions. In fact, by virtue of equality
\begin{equation*}
{\sum\limits_{k = 2}^{n} {{\sum\limits_{i = 2}^{k} {m_{i,k}}} } }  =
{\sum\limits_{k = 1}^{n - 1} {{\sum\limits_{i = k + 1}^{n} {m_{k + 1,i}}} }
}
\end{equation*}
\noindent
we obtain
\begin{equation}
\label{eq14}
{\sum\limits_{k = 1}^{n} {B(k,n)}}  = 2{\sum\limits_{k = 2}^{n}
{{\sum\limits_{i = 2}^{k} {m_{i,k}}} } }  = 2{\sum\limits_{k = 1}^{n - 1}
{{\sum\limits_{i = k + 1}^{n} {m_{k + 1,i}}} } } .
\end{equation}

We present some simple properties of the functions $A\left(
{k,n} \right)$ and $B\left( {k,n} \right)$, defined by the formula (9):
\begin{equation}
\label{eq15}
A\left( {n + 1,n + 1} \right) - B\left( {n + 1,n + 1} \right) = A\left(
{n,n} \right),
\end{equation}
\begin{equation}
\label{eq16}
A\left( {k + 1,k + 1} \right) - B\left( {k + 1,k + 1} \right) = A\left(
{k,n} \right) - B\left( {k,n} \right) + m_{2,n + 1} + ... + m_{k,n + 1} .
\end{equation}

Those properties are easily proved if we proceed from the
definitions of functions $A\left( {k,n} \right)$ and $B\left( {k,n}
\right)$.

\bigskip

\textbf{Lemma 2}. Let $a,b_{1} ,$\ldots , $b_{n} $ are real numbers with $a = 0,\,
- 1,\, - 2,...$and $a > b_{1} + ... + b_{n} .$ Then the following summation
formula holds true at $n \in {\rm N}\backslash \{1\}$
\begin{equation}
\label{eq8}
{\sum\limits_{{\mathop {m_{i,j} = 0}\limits_{(2 \le i \le j \le n)}
}}^{\infty}  {{\displaystyle\frac{{(\alpha)_{A (n,n)}}} {{m_{ij}!} }}}} {\prod\limits_{k = 1}^{n} {{\frac{{\left( {b_{k}}  \right)_{B(k,n)}
\left( {a - b_{k}}  \right)_{A(k,n) - B(k,n)}}} {{\left( {a}
\right)_{A(k,n)}}} }}} = \Gamma \left( {a - {\sum\limits_{k = 1}^{n} {b_{k}
}}}  \right){\frac{{\Gamma ^{n - 1}\left( {a} \right)}}{{{\prod\limits_{k =
1}^{n} {\Gamma \left( {a - b_{k}}  \right)}}} }}.
\end{equation}

\bigskip

Note that if we put $n = 2$ in the formula (\ref{eq8}), then

\begin{equation}
\label{eq9}
F\left( {b_{1} ,b_{2} ;a;1} \right) = {\sum\limits_{m_{22} = 0}^{\infty}
{{\frac{{\left( {b_{1}}  \right)_{m_{22}}  \left( {b_{1}}  \right)_{m_{22}}
}}{{\left( {c} \right)_{m_{22}}  m_{22} !}}}}}  = {\frac{{\Gamma \left( {a -
b_{1} - b_{2}}  \right)\Gamma \left( {a} \right)}}{{\Gamma \left( {a - b_{1}
} \right)\Gamma \left( {a - b_{2}}  \right)}}},
\end{equation}
\noindent
that is, the formula (13) is a natural generalization of the well-known
summation formula for the Gauss hypergeometric function.

\bigskip

The proof of \textbf{Lemma 2} is carried out by the method of mathematical induction.

From equality (14) it follows that the formula (13) is valid for $n = 2$.

Now we denote the left side of the formula (13) by

$$T_{n} \left( {a,b_{1} ,...,b_{n}}  \right): =
{\sum\limits_{{\mathop {m_{i,j} = 0}\limits_{(2 \le i \le j \le n)}
}}^{\infty}  {{\displaystyle\frac{{\left( {a} \right)_{A(n,n)}}} {{m_{ij}
!}}}}} {\prod\limits_{k = 1}^{n} {{\displaystyle\frac{{\left( {b_{k}}  \right)_{B(k,n)}
\left( {a - b_{k}}  \right)_{A(k,n) - B(k,n)}}} {{\left( {a}
\right)_{A(k,n)}}} }}} $$

\noindent
and considering fair equality

$$T_{n} \left( {a,b_{1} ,...,b_{n}}  \right) = \Gamma \left( {a -
{\sum\limits_{k = 1}^{n} {b_{k}}} }  \right){\frac{{\Gamma ^{n - 1}\left(
{a} \right)}}{{{\prod\limits_{k = 1}^{n} {\Gamma \left( {a - b_{k}}
\right)}}} }},$$

\noindent
we will prove that

\begin{equation}
\label{eq10}
T_{n + 1} \left( {a,b_{1} ,...,b_{n + 1}}  \right) = \Gamma \left( {a -
{\sum\limits_{k = 1}^{n + 1} {b_{k}}} }  \right){\frac{{\Gamma ^{n}\left(
{a} \right)}}{{{\prod\limits_{k = 1}^{n + 1} {\Gamma \left( {a - b_{k}}
\right)}}} }}.
\end{equation}

For this aim we will put

\[
T_{n + 1} \left( {a,b_{1} ,...,b_{n + 1}}  \right) =
{\sum\limits_{{\mathop {m_{i,j} = 0}\limits_{(2 \le i \le j \le n+1)}
}}^{\infty}   {{\frac{{\left( {a} \right)_{A(n + 1,n + 1)}
}}{{m_{ij} !}}}}} {\prod\limits_{k = 1}^{n + 1} {{\frac{{\left( {b_{k}}
\right)_{B(k,n + 1)} \left( {a - b_{k}}  \right)_{A(k,n + 1) - B(k,n + 1)}
}}{{\left( {a} \right)_{A(k,n + 1)}}} }}}
\]

\noindent
and show the validity of the recurrence relation
\begin{equation}
\label{eq11}
T_{n + 1} (a,b_{1} ,...,b_{n + 1} ) = {\prod\limits_{k = 1}^{n} {{\left[
{{\frac{{\Gamma \left( {a} \right)\Gamma \left( {a - b_{k} - b_{n + 1}}
\right)}}{{\Gamma \left( {a - b_{n + 1}}  \right)\Gamma \left( {a - b_{k}}
\right)}}}} \right]}}} \,T_{n} (a - b_{n + 1} ,b_{1} ,...,b_{n} ).
\end{equation}

This process consists of $n$ steps. A detailed look at the first step.

By virtue of the equalities
\[
{\sum\limits_{{\mathop {m_{i,j} = 0}\limits_{(2 \le i \le j \le n+1)}
}}^{\infty}   {}}  = {\sum\limits_{{\mathop {m_{i,j} = 0}\limits_{(2 \le i \le j \le n)}
}}^{\infty}   {{\sum\limits_{{\mathop {m_{i,n+1} = 0}\limits_{(2 \le i \le n+1)}
}}^{\infty}   {}}} }  = {\sum\limits_{{\mathop {m_{i,j} = 0}\limits_{(2 \le i \le j \le n)}
}}^{\infty}   {{\sum\limits_{{\mathop {m_{i,n+1} = 0}\limits_{(2 \le i \le j \le n)}
}}^{\infty}  }}\sum\limits_{m_{n+1,n+1} = 0}
^{\infty} }
\]
\noindent
and the properties of functions $A\left( {k,n} \right)$ and $B\left( {k,n}
\right)$ (see formulas (11) and (12)), the right side of equality

$$T_{n + 1} \left( {a,b_{1} ,...,b_{n + 1}}  \right) =
{\sum\limits_{{\mathop {m_{i,j} = 0}\limits_{(2 \le i \le j \le n+1)}
}}^{\infty} {{\frac{{\left( {a} \right)_{A(n + 1,n + 1)}
}}{{m_{ij} !}}}}} {\prod\limits_{k = 1}^{n + 1} {{\frac{{\left( {b_{k}}
\right)_{B(k,n + 1)} \left( {a - b_{k}}  \right)_{A(k,n + 1) - B(k,n + 1)}
}}{{\left( {a} \right)_{A(k,n + 1)}}} }}} $$

\noindent
it is easy to convert to the form
$$
\begin{array}{l}
 T_{n + 1} \left( {a,b_{1} ,...,b_{n + 1}}  \right) =
{\sum\limits_{{\mathop {m_{i,j} = 0}\limits_{(2 \le i \le j \le n)}
}}^{\infty}  {{\displaystyle\frac{{\left( {a - b_{n + 1}}  \right)_{A(n,n)}
\left( {b_{n}}  \right)_{B(n,n)}}} {{m_{ij} !}}}}}  \\
\\
 {\sum\limits_{{\mathop {m_{i,n+1} = 0}\limits_{(2 \le i  \le n)}
}}^{\infty}  {{\displaystyle\frac{{\left( {b_{n + 1}}  \right)_{m_{2,n + 1} +
... + m_{n,n + 1}}  \left( {a - b_{n}}  \right)_{A(n,n) - B(n,n) + m_{2,n +
1} + ... + m_{n,n + 1}}} } {{m_{i,n + 1} !\left( {a} \right)_{A(n,n) +
m_{2,n + 1} + ... + m_{n,n + 1}}} } }}}  \\
\\
 {\prod\limits_{k = 1}^{n - 1} {{\left[ {{\displaystyle\frac{{\left( {b_{k}}
\right)_{B(k,n) + m_{k + 1,n + 1}}  \left( {a - b_{k}}  \right)_{A(k,n) -
B(k,n) + m_{2,n + 1} + ... + m_{k,n + 1}}} } {{\left( {a} \right)_{A(k,n) +
m_{2,n + 1} + ... + m_{k + 1,n + 1}}} } }S(k,n)} \right]}}} , \\
 \end{array}
$$

\noindent
where
$$
S(k,n) = {\sum\limits_{m_{n + 1,n + 1} = 0}^{\infty}  {{\frac{{\left( {b_{n}
+ B(n,n)} \right)_{m_{n + 1,n + 1}}  \left( {b_{n + 1} + m_{2,n + 1} + ... +
m_{n,n + 1}}  \right)_{m_{n + 1,n + 1}}} } {{m_{n + 1,n + 1} !\left( {a +
A(n,n) + m_{2,n + 1} + ... + m_{n,n + 1}}  \right)_{m_{n + 1,n + 1}}} } }}
}.
$$

It is easy to notice that
$$
S(k,n) = F\left( b_{n} + B(n,n),b_{n + 1} + m_{2,n + 1} + ... + m_{n,n + 1}
;\right. $$
$$ \left. a + A(n,n) + m_{2,n + 1} + ... + m_{n,n + 1} ;1 \right).
$$

Applying now the summation formula (\ref{eq9}) to the last equality after
elementary transformations we get
$$
 T_{n + 1}^{(1)} \left( {a,b_{1} ,...,b_{n + 1}}  \right) = {\frac{{\Gamma
\left( {a - b_{n} - b_{n + 1}}  \right)\Gamma \left( {a} \right)}}{{\Gamma
\left( {a - b_{n}}  \right)\Gamma \left( {a - b_{n + 1}}
\right)}}}
$$
$$
{\sum\limits_{{\mathop {m_{i,j} = 0}\limits_{(2 \le i \le j \le n+1)}
}}^{\infty} {{\frac{{\left( {b_{n}}  \right)_{B(n,n)} \left( {a
- b_{n} - b_{n + 1}}  \right)_{A(n,n) - B(n,n)}}} {{m_{ij} !}}}}} {\sum\limits_{{\mathop {m_{i,n+1} = 0}\limits_{(2 \le i \le j \le n)}
}}^{\infty}  {{\frac{{\left( {b_{n + 1}}  \right)_{m_{2,n + 1} +
... + m_{n,n + 1}}} } {{m_{i,n + 1} !}}}}}
$$
$$
{\prod\limits_{k = 1}^{n - 1}
{{\frac{{\left( {b_{k}}  \right)_{B(k,n) + m_{k + 1,n + 1}}  \left( {a -
b_{k}}  \right)_{A(k,n) - B(k,n) + m_{2,n + 1} + ... + m_{k,n + 1}}
}}{{\left( {a} \right)_{A(k,n) + m_{2,n + 1} + ... + m_{k + 1,n + 1}}} } }}
}.
$$

For definiteness, we denote the result of the first step of the process
under consideration by $T_{n + 1}^{(\ref{eq1})} \left( {a,b_{1} ,...,b_{n + 1}}
\right)$. We continue the process of proving the recurrence relation (\ref{eq11}).
In each next step, having consistently repeated the reasoning carried out in
the first step, we get
$$
 T_{n + 1}^{(s)} \left( {a,b_{1} ,...,b_{n + 1}}  \right) = {\frac{{\Gamma
^{s}\left( {a} \right)}}{{\Gamma ^{s}\left( {a - b_{n + 1}}
\right)}}}{\prod\limits_{k = n - s + 1}^{n} {{\frac{{\Gamma \left( {a -
b_{k} - b_{n + 1}}  \right)}}{{\Gamma \left( {a - b_{k}}  \right)}}}}}
 $$
 $$
 {\sum\limits_{{\mathop {m_{i,j} = 0}\limits_{(2 \le i \le j \le n)}
}}^{\infty} {{\frac{{1}}{{m_{ij} !}}}{\prod\limits_{k = n - s +
1}^{n} {{\left[ {{\frac{{\left( {b_{k}}  \right)_{B(k,n)} \left( {a - b_{k}
- b_{n + 1}}  \right)_{A(k,n) - B(k,n)}}} {{\left( {a - b_{n + 1}}
\right)_{A(k,n)}}} }} \right]}}}} }
$$
$$
 {\sum\limits_{{\mathop {m_{i,n+1} = 0}\limits_{(2 \le i  \le n-s+1)}
}}^{\infty} {}} {\frac{{\left( {a - b_{n + 1}}  \right)_{N(n,n)}
\left( {b_{n + 1}}  \right)_{m_{2,n + 1} + ... + m_{n - s + 1,n + 1}}
}}{{m_{ij} !}}}
$$
$$
 {\prod\limits_{k = 1}^{n - s} {{\left[ {{\frac{{\left( {b_{k}}
\right)_{B(k,n) + m_{k + 1,n + 1}}  \left( {a - b_{k}}  \right)_{A(k,n) -
B(k,n) + m_{2,n + 1} + ... + m_{k,n + 1}}} } {{\left( {a} \right)_{A(k,n) +
m_{2,n + 1} + ... + m_{k + 1,n + 1}}} } }} \right]}}}
$$
\noindent
and in the last step
$$
T_{n + 1}^{(n)} \left( {a,b_{1} ,...,b_{n + 1}}  \right) = {\frac{{\Gamma
^{n}\left( {a} \right)}}{{\Gamma ^{n}\left( {a - b_{n + 1}}
\right)}}}{\prod\limits_{k = 1}^{n} {{\left[ {{\frac{{\Gamma \left( {a -
b_{n + 1} - b_{k}}  \right)}}{{\Gamma \left( {a - b_{k}}  \right)}}}}
\right]}}}
$$
$$
 {\sum\limits_{{\mathop {m_{i,j} = 0}\limits_{(2 \le i \le j \le n)}
}}^{\infty}  {{\frac{{\left( {a - b_{n + 1}}  \right)_{A(n,n)}
}}{{m_{ij} !}}}}} {\prod\limits_{k = 1}^{n} {{\left[ {{\frac{{\left( {b_{k}
} \right)_{B(k,n)} \left( {a - b_{n + 1} - b_{k}}  \right)_{A(k,n) - B(k,n)}
}}{{\left( {a - b_{n + 1}}  \right)_{A(k,n)}}} }} \right]}}} , \\
$$
\noindent
that is
$$
T_{n + 1}^{(n)} \left( {a,b_{1} ,...,b_{n + 1}}  \right) = {\frac{{\Gamma
^{n}\left( {a} \right)}}{{\Gamma ^{n}\left( {a - b_{n + 1}}
\right)}}}{\prod\limits_{k = 1}^{n} {{\left[ {{\frac{{\Gamma \left( {a -
b_{n + 1} - b_{k}}  \right)}}{{\Gamma \left( {a - b_{k}}  \right)}}}}
\right]}}} T_{n} \left( {a - b_{n + 1} ,b_{1} ,...,b_{n}}  \right).
$$

Thus, the validity of the ratio (\ref{eq11}) is established. By the induction
hypothesis, from the (\ref{eq11}) follows the equality
$$
T_{n} \left( {a - b_{n + 1} ,b_{1} ,...,b_{n}}  \right) = \Gamma \left( {a -
b_{n + 1} - {\sum\limits_{k = 1}^{n} {b_{k}}} }  \right){\frac{{\Gamma ^{n -
1}\left( {a - b_{n + 1}}  \right)}}{{{\prod\limits_{k = 1}^{n} {\Gamma
\left( {a - b_{n + 1} - b_{k}}  \right)}}} }}.
$$

Substituting the last expression in (\ref{eq11}) we get equality (\ref{eq8}). Q.E.D.

\bigskip

\textbf{Lemma 3}. The following equality
$$
{\mathop{\lim} \limits_{\mathop{z_{k} \to 0,}\limits_{k =1,...,n}}} z_{1}^{ - b_{1}}  ...z_{n}^{ - b_{n}}  F_{A}^{(n)} \left(
{a,b_{1} ,...,b_{n} ;c_{1} ,...,c_{n} ;1 - {\frac{{1}}{{z_{1}}} },...,1 -
{\frac{{1}}{{z_{n}}} }} \right)
$$
\begin{equation}
\label{eq12}
= {\frac{{1}}{{\Gamma (a)}}}\Gamma \left( {a - {\sum\limits_{k = 1}^{n}
{b_{k}}} }  \right){\prod\limits_{k = 1}^{n} {{\frac{{\Gamma \left( {c_{k}}
\right)}}{{\Gamma \left( {c_{k} - b_{k}}  \right)}}}}}
\end{equation}
\noindent
is valid.

\textbf{Proof}. By virtue of the decomposition formula (8) we obtain
\[
F_{A}^{(n)} \left( {a,b_{1} ,...,b_{n} ;c_{1} ,...,c_{n} ;1 -
{\frac{{1}}{{z_{1}}} },...,1 - {\frac{{1}}{{z_{n}}} }} \right) =
{\sum\limits_{{\mathop {m_{i,j} = 0}\limits_{(2 \le i \le j \le n)}
}}^{\infty}  {{\frac{{(a)_{A(n,n)}}} {{m_{ij}!} }}}}
\]
\begin{equation}
\label{eq13}
 \cdot {\prod\limits_{k = 1}^{n}\left[ {{\frac{{(b_{k} )_{B(k,n)}}} {{(c_{k}
)_{B(k,n)}}} }\left( {1 - {\frac{{1}}{{z_{k}}} }} \right)^{B(k,n)}F\left( {a
+ A(k,n),b_{k} + B(k,n);c_{k} + B(k,n);1 - {\frac{{1}}{{z_{k}}} }} \right)}\right]
}.
\end{equation}

Applying now the familiar autotransformation formula
\[
F\left( {a,b;c;x} \right) = (1 - x)^{ - b}F\left( {c - a,b;c;{\frac{{x}}{{x
- 1}}}} \right)
\]

\noindent
to each hypergeometric function included in the sum (\ref{eq13}), we get
\[
F_{A}^{(n)} \left( {a,b_{1} ,...,b_{n} ;c_{1} ,...,c_{n} ;1 -
{\frac{{1}}{{z_{1}}} },...,1 - {\frac{{1}}{{z_{n}}} }} \right) =
z_{1}^{b_{1}}  ...z_{n}^{b_{n}}  {\sum\limits_{{\mathop {m_{i,j} =
0}\limits_{(2 \le i \le j \le n)}}} ^{\infty}  {{\frac{{(a)_{A(n,n)}
}}{{m_{ij}!} }}}}
\]

\[
 \cdot \prod\limits_{k = 1}^{n} \left[{{\frac{{(b_{k} )_{B(k,n)}}} {{(c_{k}
)_{B(k,n)}}} }\left( {z_{k} - 1} \right)^{B(k,n)}F {\left(
{{\begin{array}{*{20}c}
 {c_{k} - a + B(k,n)
- A(k,n),b_{k} + B(k,n)} ; \hfill \\
 {c_{k} + B(k,n)};  \hfill \\
\end{array}} 1 - z_{k}}  \right)}}\right].
\]

Next, we calculate the limit
$$
{\mathop{\lim} \limits_{\mathop{z_{k} \to 0,}\limits_{k =1,...,n}}} z_{1}^{ - b_{1}}  ...z_{n}^{ - b_{n}}  F_{A}^{(n)} \left(
{a,b_{1} ,...,b_{n} ;c_{1} ,...,c_{n} ;1 - {\frac{{1}}{{z_{1}}} },...,1 -
{\frac{{1}}{{z_{n}}} }} \right)
$$
\noindent
and the resulting expression we apply the summation formula (\ref{eq9}), with the
result that we obtain the equality (\ref{eq12}). Q.E.D.

 \end{document}